\documentclass[11pt]{article}
\usepackage[paper=a4paper, margin=28mm]{geometry}
\usepackage{amssymb,amsmath,amsthm}
\usepackage{amsfonts,amsmath,amssymb,graphicx}
\usepackage{tikz}
\usepackage{float}
\usepackage[toc,page]{appendix}
\usepackage{caption}
\usepackage{subcaption}
\usepackage{algorithm}
\usepackage{algpseudocode}
\usepackage{parskip}
\makeatletter
\def\thm@space@setup{%
\thm@preskip=\parskip \thm@postskip=0pt }

\renewenvironment{proof}[1][\proofname]{\par
  \pushQED{\qed}%
  \normalfont \topsep6\p@\@plus6\p@\relax
  \trivlist
  \item[\hskip\labelsep
        \itshape
    #1\@addpunct{.}]\ignorespaces\unskip
}{%
  \popQED\endtrivlist\@endpefalse
}
\makeatother

\newcommand{\Erdos}{Erd\H{o}s }

\makeatletter
\renewcommand{\maketitle}
{\begin{center}
\LARGE \@title\\
\vspace{1em} \@author \vspace{1em}
\end{center}
}
\makeatother

\newtheorem{thm}{Theorem}[section]
\newtheorem{lemma}[thm]{Lemma}
\newtheorem{prop}[thm]{Proposition}
\newtheorem{cor}[thm]{Corollary}
\theoremstyle{definition}
\newtheorem{defn}[thm]{Definition}
\newtheorem{example}[thm]{Example}
\begin{document}
\title{On sum-free subsets of abelian groups}
\author{Renato Cordeiro de Amorim}\footnote{School of Computer Science and Electronic Engineering, University of Essex, CO4 3SQ, United Kingdom. Contact e-mail: r.amorim@essex.ac.uk}
\maketitle


\begin{abstract}
In this paper we discuss some of the key properties of sum-free subsets of abelian groups. Our discussion has been designed with a broader readership in mind, and is hence not overly technical. We consider answers to questions like: how many sum-free subsets are there in a given abelian group $G$? what are its sum-free subsets of maximum cardinality? what is the maximum cardinality of these sum-free subsets? what does a typical sum-free subset of $G$ looks like? among others.
%
%
\\\textbf{Keywords}: sum-free sets; abelian groups.

\end{abstract}

\section{Introduction}
\label{chap:introduction}

Here, we present the key progress made on sum-free subsets of abelian groups. Our main aim is to convey this information without being overly technical, in order to make the topic accessible to a broader readership. 
Let us begin with a definition.

\begin{defn}
Let $G$ be an abelian group. A subset $S$ of $G$ is sum-free if for all $x,y \in S$ we have that $x+y \not \in S$.
\end{defn}

In the case of $\mathbb{Z}$, under the usual addition, one can easily come up with simple examples. For instance, $\{1,4\}$, $\{2,7,8\}$, and any set composed solely of odd numbers. Schur \cite{schur1916} in 1916 was the first to start work that would eventually lead to the study of sum-free sets. He did so by proving that for a sufficiently high $n \in \mathbb{Z}^+$ every finite colouring of the integers in the interval $[1,n]$ contains a monochromatic triple $(x,y,z)$, usually referred as a Schur triple, such that $x+y = z$ \cite{cameron1987portrait}. Cameron and \Erdos \cite{cameron1990number} also make a considerable contribution to this topic, focusing in particular on sum-free subsets of positive integers (for details see Section \ref{chapter:counting_sum_free}).

The literature on sum-free sets presents very interesting questions. For instance, Alon et al. \cite{alon2014counting} and Green and Ruzsa \cite{green2005sum} pose the following: how many sum-free subsets are there in $G$? what are its sum-free subsets of maximum cardinality? what is the maximum cardinality of these sum-free subsets? what does a typical sum-free subset of $G$ looks like? Here, we address these questions among others.

We organised our paper as follows. Section \ref{chap:Identifying_sum_free_sets} presents basic methods, computational and algebraic, capable of identifying sum-free sets. There we also discuss periodic sum-free sets, presenting methods to generate infinite sum-free subsets of positive integers. Section \ref{chapter:counting_sum_free} considers the bounds on the quantity of sum-free sets. It first discusses these bounds on subsets of positive integers, and then provides a more general account regarding any finite abelian group. Section \ref{chap:maximal} deals with maximal sum-free sets, that is, a broader class of sum-free sets that includes sum-free sets of maximum cardinality. Section \ref{chapter:conclusion} reviews our main findings and presents our conclusions.

\section{Identifying sum-free sets}
\label{chap:Identifying_sum_free_sets}

In this section we are particularly interested in the following question: given an abelian group $G$, finite or not, are there methods to identify its sum-free subsets? before we can explore potential answers it is helpful to present the bounds that follow. In the case of sum-free integer subsets of $[1,n]$, \Erdos \cite{erdijs1965extremal} shows this interval contains a sum-free set of cardinality $\frac{1}{3}n$. Alon and Kleitman \cite{alon1990sum} go even further by devising the more general theorem below.

\begin{thm}
\label{thm:existence_of_sum_free_ndiv_2_7}
Let $G$ be a finite abelian group. For every $B\subset G\setminus \{e\}$ there exists a sum-free subset of $B$ with cardinality greater than $\frac{2}{7}|B|$.
\end{thm}

Alon and Kleitman prove the above by first observing that the sets $A_1=\{x \in \mathbb{Z}_n \mid \frac{1}{3}n < x \leq \frac{2}{3}n\}$ and $A_2=\{x \in \mathbb{Z}_n \mid \frac{1}{6}n < x \leq \frac{1}{3}n \text{ or } \frac{2}{3}n <x\leq n\}$ are sum-free subsets of $\mathbb{Z}_n$. Given $\mathbb{Z}_n$ is cyclic, it has a subgroup $d\mathbb{Z}_n = \{0, d, 2d, \ldots, n-d\}$ for every $d$ dividing $n$, with order $\frac{n}{d}$. They then provide a case analysis computing $\frac{|d\mathbb{Z}_n \cap A_i|}{|d\mathbb{Z}_n|}$ for all divisors $d$ of $n$ and $i\in \{1,2\}$, finding that 
\begin{align*}
    \frac{4}{7}\frac{|d\mathbb{Z}_n \cap A_1|}{|d\mathbb{Z}_n|}+\frac{3}{7}\frac{|d\mathbb{Z}_n \cap A_2|}{|d\mathbb{Z}_n|} \geq \frac{2}{7}.
\end{align*}

The above eventually leads to the bound in Theorem \ref{thm:existence_of_sum_free_ndiv_2_7}. Alon and Kleitman then make use of a theorem by Rhemtulla and Street \cite{rhemtulla1971maximal} regarding the maximum cardinality of sum-free subsets of elementary abelian groups to show that $\frac{2}{7}$ is in fact optimal. 

Theorem \ref{thm:existence_of_sum_free_ndiv_2_7} provides a lower bound for the cardinality of the largest sum-free subset of $G$. This can be helpful if one is trying to identify all of its sum-free subsets. That is, if no sum-free set of cardinality greater than $\frac{2}{7}|B|$ has been identified, then the list is not complete.

The rest of this section discusses methods to help with the identification of sum-free sets. We explore three simple algorithms (for completeness) capable of generating sum-free sets, periodic sum-free sets, and other useful algebraic methods.

\subsection{Computational approaches}
\label{sec:computational_approaches_find_sum_free}

Cameron \cite{cameron_blog} describes an algorithm to identify sum-free subsets of positive integers. As he puts: when examining a number $a$, if this is not the sum of two elements of a set $S$, then choose uniformly at random whether to put $a$ in $S$. This informal description does not clarify how $a$ is selected, or what criteria should be used to stop the algorithm. Although these points may seem negative, his approach does make the algorithm easily adaptable to particular cases. 

We formalise the steps above in Algorithm \ref{alg:prob_sum_free}. We decided to set the first value $a\in \mathbb{Z}^+$ to be user-defined (it could be chosen at random). We also introduced a condition regarding the cardinality of the sum-free set to be identified, to ensure an eventual completion. This naive algorithm is capable of identifying a sequence of nested sum-free subsets of $\mathbb{Z}^+$. If applied to a finite group it could be easily modified to identify maximal sum-free sets (for details on the latter, see Section \ref{chap:maximal}), this is what we meant with adaptable. Unfortunately, this algorithm will not identify all sum-free sets. For instance, it is unable to identify two sum-free sets with an empty intersection.

\begin{algorithm}
\caption{Cameron's algorithm to identify sum-free sets}
\label{alg:prob_sum_free}
\begin{algorithmic}
\Require $a, b \in \mathbb{Z}^+$, where $a$ is the initial element to be considered and $b$ is the cardinality of the sum-free subset of $\mathbb{Z}^+$ to be identified.
\Ensure $S$ is a sum-free subset of $\mathbb{Z}^+$ such that $a \in S$.
\State $S \gets \emptyset$
\While{$|S| < b$}
    \State Identify an element $s \in \mathbb{Z}^+$ uniformly at random.
    \State $A \gets S \cup \{s\}$.
    \State $v \gets $ rand($\{1,2\}$).
    \If{$A \cap (A +A)=\emptyset$ and $v=1$}
        \State $S \gets A$.
    \EndIf
\EndWhile
\end{algorithmic}
\end{algorithm}


With the above in mind one may wonder if it is possible to have an algorithm capable of identifying all sum-free sets. The obvious approach would be to test each of the possible subsets of $G$ and list those that are sum-free. However, this is rather computationally expensive running at $\mathcal{O}(2^{|G|-1})$. Algorithm \ref{alg:all_sum_free}, where we use a sequence of binary numbers to identify each subset, formalises the steps. Notice this algorithm can be made parallel by supplying $k$ computers with different sets of binary numbers.

\begin{algorithm}
\caption{Identifying all sum-free sets}
\label{alg:all_sum_free}
\begin{algorithmic}
\Require An additive non-trivial finite abelian group $G$, with identity $e$.
\Ensure $S$ is a set containing all sum-free subsets of $G$.
\State $S \gets \emptyset$
\State Put each element of $G\setminus \{e\}$ in a vector $\textbf{v}$.
\For{$i=1$ to $2^{|G|-1}$}
    \State Create a vector $\textbf{b}$ containing the binary value of $i$ spread over its $|G|-1$ components.
    \State Create a set $A$ such that the element $\textbf{v}[j] \in A$ if \textbf{b}$[j]=1$, for $j=1\ldots, |G|-1$.
    \If{$A \cap (A +A)=\emptyset$}
        \State Put $A$ in $S$.
    \EndIf
\EndFor
\end{algorithmic}
\end{algorithm}

Taking into account the prohibitively high exponential bound above Kolountzakis \cite{kolountzakis1994selection} provides a partial solution that applies to integer sum-free sets. We have that every set $A=\{a_1, \ldots, a_n\}$ of integers contains a sum-free subset $S$ with $|S|>\frac{2}{7}n$ (see Theorem \ref{thm:existence_of_sum_free_ndiv_2_7}). Kolountzakis proposes an algorithm capable of identifying $S$ in polynomial time. Notice that $S$ is a large sum-free set that will contain more sum-free sets if $|S|>1$ (any subset of a sum-free set is sum-free) but it may not contain all. Hence, our use of the word partial. Kolountzakis' algorithm is based on the following theorem.

\begin{thm}
Let $p=3k+2$ be a prime and $w=\sum_{x \in \mathbb{Z}_p^{\times}} w(x)$, where $w(x)$ is non-negative and defined on $\mathbb{Z}_p^{\times}$. Then, there exists a sum-free set $S\in \mathbb{Z}_p^{\times}$ such that $\sum_{x \in S} w(x) > \frac{1}{3}w$.
\end{thm}

In the above $\mathbb{Z}_p^{\times}$ refers to the multiplicative group over $\mathbb{Z}_p$. Kolountzakis proves the above by first noticing that $B=\{k+1, \ldots, 2k+1\}$ is sum-free in $\mathbb{Z}_p$ and $|B| > \frac{p-1}{3}$. Kolountzakis then builds on $B$ by identifying an element $t \in \mathbb{Z}_p$ for which $\sum_{t \cdot x \in B}w(x)>\frac{w}{3}$. Finally, he defines the sum-free set $S=t^{-1}B$ proving $S$ is sum-free and has the required cardinality. Algorithm \ref{alg:large_sum_free} presents the steps to identify $S$, which defines $w(x) = |\{a \in A \mid a= x \text{ (mod } p\text{)}\}|$ with $A$ being a set of integers.

Kolountzakis makes some considerations regarding the computational cost of his algorithm. He notes that the number of prime factors of an integer $x$ is at most $\log_2 x$. Hence, the number of primes in the factorisation of any $x \in A$ is at most $l=\sum_{j=1}^n \log_2 a_j$. He then applies the Prime Number Theorem of arithmetic progressions to state that there is a prime $p=3k+2 \leq 3l \log l$, which does not divide any $a \in A$, concluding the algorithm can be carried out in polynomial time.

\begin{algorithm}
\caption{Identifying a large sum-free set (Kolountzakis)}
\label{alg:large_sum_free}
\begin{algorithmic}
\Require A set of integers $A=\{a_1, \ldots, a_n\}$.
\Ensure $S$ is a sum-free subset of $A$.
\State Identify a prime $p=3k+2$ with $p\leq 3l \log l$ such that it divides no $a \in A$.
\State Compute $w(x) = |\{a \in A \mid a=x \text{ mod } p\}|$ for all $x \in \mathbb{Z}_p^{\times}$.
\State Identify a $t \in \mathbb{Z}_p^{\times}$ for which $\sum_{t\cdot x \in B}w(x)>\frac{n}{3}$.
\State Construct the set $S^{\prime} = t^{-1}B$
\State Construct the set $S = \{x \in A \mid x \text{ mod } p \in S^{\prime}\}$.
\end{algorithmic}
\end{algorithm}

\subsection{Periodic sum-free sets}
\label{sec:periodic_sum_free}

The set of odd integers is probably the most well-known example of a periodic sum-free set. There are, of course, other periodic sum-free subsets of integers.

\begin{lemma}
\label{lemma:sum_free_xequivm_modn}
Let $S=\{x \in \mathbb{Z} \mid x\equiv m \text{ (mod n)}\}$ with $0<m<n$, $S$ is sum-free.
\end{lemma}
\begin{proof}
Let $x,y \in S$. Then, $x=kn+m$ and $y=qn+m$ for some $k,q \in \mathbb{Z}$. Let us assume, for a contradiction, that $x+y \in S$. This would mean $(kn+m)+(qn+m)-m=tn$ for some $t \in \mathbb{Z}$, which can be simplified to $m=n(t-k-q)$. However, the latter is a contradiction as we have $0<m<n$.
\end{proof}

Tran \cite{tran2018structure} states that the sets $\{x \in \mathbb{Z} \mid x\equiv 1,4 \text{ (mod 5)}\}$ and $\{x \in \mathbb{Z} \mid x\equiv 2,3 \text{ (mod 5)}\}$ are sum-free. This is indeed the case, and we can devise another more general result.

\begin{lemma}
Let $A$ be a sum-free subset of $\mathbb{Z}_n$ and $S=\{x \in \mathbb{Z} \mid \exists a \in A \text{ s.t. } x\equiv a \text{ (mod n)}\}$, $S$ is sum-free.
\end{lemma}
\begin{proof}
Let $x,y \in S$. Then, $x=kn+a_i$ and $y = qn+a_j$ for some $k,q \in \mathbb{Z}$ and $a_i,a_j \in A$. We have that $x+y=(k+q)n+a_i+a_j$. Given $a_i \oplus_n a_j \not \in A$, we have that $x+y \not \in S$.
\end{proof}

\subsection{Other algebraic approaches}
\label{sec:algebraic_approached}

In this section we explore other algebraic approaches to identify sum-free subsets of $G$. This tasks becomes considerably easier if one already knows one such, preferably large, sum-free set.

\begin{prop}
 \label{prop:subset_of_sumfree_is_sumfree}
 Let $S$ be a sum-free subset of $G$, then any $A \subseteq S$ is sum-free.
\end{prop}
\begin{proof}
$A$ can only contain elements that are in $S$. Hence, $A$ is sum-free.
\end{proof}

\begin{prop}
\label{prop:_intersection_sumfree_is_sumfree}
Let $A$ and $B$ be sum-free subsets of an abelian group $G$. The intersection of $A$ and $B$ is sum-free.
\end{prop}
\begin{proof}
We have that $A\cap B \subseteq A, B$. Hence, Proposition \ref{prop:subset_of_sumfree_is_sumfree} tell us $A \cap B$ is sum-free.
\end{proof}

Proposition \ref{prop:_intersection_sumfree_is_sumfree} has an interesting consequence, first pointed out by Cameron and \Erdos \cite{cameron1990number} for the case of integer sum-free subsets of $[1,n]$, which we show in a more general form below.

\begin{cor}
\label{cor:sum_free_closed_under_intersection}
Let $\mathcal{F}$ be a family of sum-free subsets of an abelian group $G$. $\mathcal{F}$ is closed under intersection.
\end{cor}

If there is no known sum-free subset for the group in question but this group is finite and a subgroup $H$ is known, then one can use the lemma below to identify sum-free sets.

\begin{lemma}
\label{lemma:some_coset_is_sum_free}
Let $H$ be a proper subgroup of a finite group $G$. Given an element $g \in G\setminus H$, the coset $H+g$ is sum-free.
\end{lemma}
\begin{proof}
Let us assume that $H+g$ is not sum-free. Then, $h_1+g = (h_2+g)+(h_3+g)$ for some $h_1, h_2, h_3 \in H$ and $g \in G\setminus H$. This implies $g=(-h_2)+h_1+(-h_3)$, and by consequence $g \in H$. However, this contradicts the hypothesis. Hence, $H+g$ is sum-free.
\end{proof}

\section{Counting sum-free sets}
\label{chapter:counting_sum_free}
We begin this section by discussing the foundational work of Cameron and \Erdos \cite{cameron1990number}, further restricting the definition of sum-free sets to focus solely on positive integers (see Sections \ref{sec:ce_conjecture}, \ref{sec:advances_ce_conjecture}, and \ref{sec:proof_ce_conjecture}). Section \ref{sec:counting_sum_free_in_general} generalises our discussion to any abelian group.

\subsection{The Cameron and \Erdos\kern-0.65ex' conjecture}
\label{sec:ce_conjecture}

Let $S$ be a sum-free set of integers, such that $S\subseteq [1,n]$ with largest element $k$. In their seminal work Cameron and \Erdos \cite{cameron1990number} start their analysis of sum-free sets like $S$ with a straightforward remark, leading to interesting consequences. For each integer $i<k$, $S$ contains at most one component of the pair $(i, k-i)$. Hence, $|S| \leq \lceil \frac{1}{2}k \rceil$ and by consequence $|S| \leq \lceil \frac{1}{2}k\rceil \leq \lceil \frac{1}{2}n\rceil$. 

The above is rather interesting because it easily gives us an upper bound for $|S|$. Cameron and \Erdos \cite{cameron1990number} carry on stating that, in the case of positive integers, the only sum-free sets of cardinality $\lceil \frac{1}{2}n \rceil$ are:
\begin{enumerate}
    \item The odd numbers in the interval $[1,n]$;
    \item If $n$ is odd, $\left[ \frac{1}{2}(n+1),n\right]$;
    \item If $n$ is even, $\left[\frac{1}{2}n,n-1\right]$ and $\left[\frac{1}{2}n+1,n\right]$.
\end{enumerate}


Any set composed solely of odd numbers is sum-free, and it is straightforward to see that any of the sets above has the required cardinality of $\lceil \frac{1}{2}n\rceil$. Freiman \cite{freiman1991structure} proves the above rules hold for $n\geq 24$. However, using Algorithm \ref{alg:all_sum_free} we have found the following exceptions for $n<24$: $n=4$, $\{1,4\}$; $n=6$, $\{2,5,6\}$ and $\{1,4,6\}$; $n=8$, $\{2,3,7,8\}$. These are the only exceptions for the rules above.

Clearly, the existence of the exceptions above does not invalidate the stated upper bound for $|S|$. Let $f(n)$ be the number of sum-free sets in $[1,n]$. On this regard, Cameron and \Erdos state that
\begin{equation*}
    \frac{f(n)}{2^{\frac{1}{2}n}}>1.
\end{equation*}

Other authors, such as Calkin \cite{calkin1990number}, reached the similar conclusion that $f(n) = 2^{\left(\frac{1}{2}+o(1)\right)n}$ but via a completely different path (for details, see Section \ref{sec:advances_ce_conjecture}). Cameron and \Erdos \cite{cameron1990number} go further and conjecture that $f(n)/2^{\frac{1}{2}n}$ is in fact bounded. They believe this ratio tends to the limits of approximately $6.8$ and $6.0$ depending on whether $n$ tends to infinity through odd or even numbers. They then introduce the following theorem.

\begin{thm}
There is an absolute constant $c$ such that the number of sum-free sets of $[1,n]$ whose least element is greater than $\frac{n}{3}$ does not exceed $c\times 2^{\frac{1}{2}n}$.
\end{thm}

The complete proof for the above can be found in \cite{cameron1990number}. Generally speaking, it begins by fixing $k \leq n/3$ and counting the sum-free subsets of $[1,2n-1]$ with no element lower than $n-k$. This proof shows that such sum-free sets fall within three categories, and that if one were to sum the quantity of all possible such sets the answer would be at most $c \times 2^n$. Given $c$ is a constant, Cameron and \Erdos conjecture that $f(n)$ is $\mathcal{O}(2^{\frac{1}{2}n})$.

\subsection{Advances on the Cameron and \Erdos\kern-0.65ex' conjecture}
\label{sec:advances_ce_conjecture}

The Cameron and \Erdos\kern-0.65ex' conjecture that $f(n)$ is $\mathcal{O}(2^{\frac{1}{2}n})$ was rather difficult to prove in full. An attempt by Freiman \cite{freiman1991structure} addressed a weaker version of this conjecture. Freiman showed that the number of sum-free sets $S \subset [1,n]$ for which $|S| \geq \frac{5}{12}k+2$ has the bound $\mathcal{O}(2^{\frac{1}{2}n})$, where $k$ is the largest element in the set. This is no quite the same as showing this bound applies to $f(n)$, that is, the number of all sum-free sets in $[1,n]$ but it is certainly a good start. His approach was to first prove the following.
\begin{thm}
\label{thm:Freiman}
If $S$ is a sum-free set of positive integers such that $|S|\geq \frac{5}{12}k+2$, where $k=\max(S)$, then either: (i) $S$ is composed solely of odd numbers, or (ii) $S$ contains both odd and even numbers, $\textup{min}(S)\geq |S|$, and $|S\cap [1,\frac{1}{2}k]| \leq \frac{k-2|S|+3}{4}$.
\end{thm}

Freiman begins his proof defining that a set $S$ is difference-free if $S \cap (S-S)=\emptyset$, where $S-S=\{a-b \mid a,b \in S \}$. He then briefly shows that sum-free and difference-free are equivalent properties of a set. Notice that $a \in S-S$ implies that there are $b,c \in S$ such that $a=b-c$. This also means that if $a$ is positive, then $-a \in S-S$ as $c-b \in S-S$. Hence, $S-S$ contains positive differences $(S-S)_+$, negative differences $(S-S)_-$, and zero. We have that $|(S-S)_+| = |(S-S)_-|$, thus $|S-S|=2|(S-S)_+|+1$. Freiman also states that $S$ and $(S-S)_+$ are subsets of $[1,k]$. This is true for $S$ by definition. 

Freiman continues by stating that if $S$ is difference-free then $|S|+|(S-S)_+| \leq k$. This is indeed the case because it could only be false if $S$ and $(S-S)_+$ shared elements, however, the definition of difference-free tells us they do not. Freiman goes to prove Theorem \ref{thm:Freiman} by analysing the different possible greatest common divisors over the elements of $S$ subtracted by min$(S)$, and using the above inequality.

Calkin \cite{calkin1990number} provided another step in proving Cameron and \Erdos\kern-0.65ex' conjecture with the theorem below.
\begin{thm}
\label{thm:calkin}
The number of sum-free sets of integers contained in $[1,n]$ is $\mathcal{O}(2^{(\frac{1}{2}+\epsilon )n})$ for every $\epsilon>0$.
\end{thm}

Again, this is not quite the same as the original conjecture that $f(n)$ is $\mathcal{O}(2^{\frac{1}{2}n})$ but its certainly close. Calkin proves the above using a theorem guaranteeing the existence of a function $g(n)$ such that every integer subset of $[1,n]$ of size $g(n)$ contains an arithmetic progression of a certain length (for details, see \cite{szemeredi1975sets}), and introducing other three lemmas. The general idea behind Calkin's proof is that it is possible to find a particular arithmetic progression shared by all integer sum-free sets of a particular size. Then, one can count how many possible sum-free sets in $[1,n]$ contain this particular arithmetic progression.

Alon \cite{alon1991independent} also made constructive steps towards proving the Cameron and \Erdos conjecture, with the following theorem.
\begin{thm}
\label{thm_alon_independent_sets_on_graphs}
Any (finite, undirected, and simple) $k$-regular graph on $n$ vertices has no more than $2^{(\frac{1}{2}+\mathcal{O}(k^{-0.1}))n}$ independent sets.
\end{thm}
The theorem above is essentially an upper bound on the number of non-adjacent vertices in a graph. Alon proves Theorem \ref{thm_alon_independent_sets_on_graphs} by showing that the upper bound for the probability of a certain event happening, on a set of mutually independent events, is exponentially small. Combining this and random two-colourings of a graph (which are mutually independent), he shows that certain colourings would not take place. Hence, the upper bound in Theorem \ref{thm_alon_independent_sets_on_graphs}. By applying Theorem \ref{thm_alon_independent_sets_on_graphs} on Cayley graphs Alon reaches the following corollary.

\begin{cor}
\label{cor:alon_sum_free_sets_bound}
Let $G$ be a group and $A\subseteq G$. Then, the number of $A$-free subsets of $G$ does not exceed $2^{(\frac{1}{2}+\mathcal{O}(|A|^{-0.1}))|G|}$.
\end{cor}

Alon defines an $A$-free set as follows. Let $S,A \subseteq G$ with $G$ finite, $S$ is $A$-free if $(S+A) \cap S=\emptyset$. That is, if there is no $a \in A$ and $s_1, s_2 \in S$ such that $s_1+a = s_2$. The concept of sum-free is just a special case of the above, in which $S=A$.

There are two important things to notice in Corollary \ref{cor:alon_sum_free_sets_bound}. First, it implies Theorem \ref{thm:calkin}, Alon presents the latter as a further corollary. Second, this corollary makes no mention of our usual $[1,n]$ interval but instead talks of a group, thus it is more general. We decided to leave it in this section rather than Section \ref{sec:counting_sum_free_in_general} (which discusses sum-free subsets of abelian groups) because of its link to the Cameron and \Erdos\kern-0.65ex' conjecture. 

\subsection{Proving the Cameron and \Erdos\kern-0.65ex' conjecture}
\label{sec:proof_ce_conjecture}

The Cameron and \Erdos\kern-0.75ex' conjecture was finally proven, independently, by Green \cite{green2004cameron} and Sapozhenko \cite{sapozhenko2003cameron}. Green did so by proving the following.
\begin{thm}
\label{thm:green_ce_conjecture}
The number of sum-free sets of integers in $[1, n]$ is asymptotically $c(n)2^{\frac{1}{2}n}$, where $c(n)$ takes two different constant values depending on the parity of $n$.
\end{thm}

Green first identifies a family $\mathcal{F}$ of subsets of $[1,n]$ with three properties:
\begin{enumerate}
    \item Each $S \in \mathcal{F}$ is almost sum-free. That is, that the number of 3-tuples $(a,b,c)$ with $a,b,c \in S$ and $a+b=c$ is $o(n^2)$. 
    \item $|\mathcal{F}| = 2^{o(n)}$. That is,  $\mathcal{F}$ does not contain too many sets.
    \item Every sum-free set of integers in $[1,n]$ is a subset of a member of $\mathcal{F}$.
\end{enumerate}

Let $S$ be a sum-free set of integers in $[1,n]$, property (3) of $\mathcal{F}$ tells us that $S \subseteq S^{\prime}$ for some $S^{\prime} \in \mathcal{F}$. Property (2) ensures that the number of sets in $\mathcal{F}$ is only $2^{o(n)}$, so the number of sum-free sets $S$ in $[1,n]$ for which $|S^{\prime}| \leq \left(\frac{1}{2} - \frac{1}{120}\right)n$ is $o(2^{\frac{1}{2}n})$. By supposing $|S^{\prime}| \geq \left(\frac{1}{2} - \frac{1}{120}\right)n$ Green is able to determine that with an $\epsilon$ that is $o(n)$, $S$ either belongs to an interval of a particular length (with an error of $32\epsilon^{\frac{1}{8}}n$ elements), or it has $54\epsilon^{\frac{1}{8}}n$ even elements. Green then shows that with $o(2^{\frac{1}{2}n})$ exceptions, all sum-free $S$ in $[1,n]$ fall under two categories:
\begin{enumerate}
    \item $S$ consists entirely of odd numbers.
    \item $S$ is contained in the interval $[\lceil \frac{n+1}{3}\rceil, n]$.
\end{enumerate}

Green then applies a result by Cameron and \Erdos \cite{cameron1990number} which estimates the number of sum-free sets in the interval $[\lceil \frac{n+1}{3}\rceil, n]$, proving Theorem \ref{thm:green_ce_conjecture} and by consequence the conjecture by Cameron and \Erdos.

The original proof of Cameron and \Erdos\kern-0.65ex' conjecture by Sapozhenko has been published in Russian \cite{sapozhenko2003cameron}, with a later version in English \cite{sapozhenko2008cameron}. In the latter, Sapozhenko reaches the same conclusion of Theorem \ref{thm:green_ce_conjecture}. To do so, he defined $f(t,n)$ to be the number of sum-free subsets in the interval $[t,n]$, with shorthand $f(n)$ for $t=1$ (matching our previous notation), and $f^1(n)$ to be the number of sum-free subsets of odd numbers also in $[1,n]$. He then proved the following theorem using a purely combinatorial approach.

\begin{thm}
\label{thm:Sapozhenko_ce_conjecture}
$f(n) \approx f(\frac{1}{3}n,n)+f^1(n)$.
\end{thm}

We have that $f^1(n)$ relates to half of all subsets in $[1,n]$, so $f^1(n)=2^{\lceil \frac{1}{2}n \rceil}$. The above implies that almost all sum-free sets in $[1,n]$ are either composed entirely of odd numbers or are contained in the interval $[\frac{1}{3}n,n]$. This is a rather similar finding to that of Green \cite{green2004cameron}. Sapozhenko points out that Theorem \ref{thm:Sapozhenko_ce_conjecture} implies Theorem \ref{thm:green_ce_conjecture}. To prove the former Sapozhenko first makes a couple of important definitions. Let $\mathcal{A}$ and $\mathcal{B}$ be families of subsets of $[1,n]$. The family $\mathcal{B}$ covers $\mathcal{A}$ if for any $A \in \mathcal{A}$ there exists a $B \in \mathcal{B}$ such that $A \subseteq B$. 

Sapozhenko defines some conditions that must hold for $\mathcal{B}$ to be called correct. He then defines $\mathcal{B}$ to be an almost correct system of containers of $\mathcal{A}$, if $\mathcal{B}$ is correct for some subfamily $\mathcal{A}^{\prime} \subseteq \mathcal{A}$ such that $|\mathcal{A}\setminus A^{\prime}| = o(2^{\frac{1}{2}n})$. With these definitions at hand we can now go back to Theorem \ref{thm:Sapozhenko_ce_conjecture}.

Let $S(n)$ and $S(t,n)$ represent the families of all sum-free sets in the intervals $[1,n]$ and $[t,n]$, respectively, and $S^1(n)$ represent the family of all sum-free sets composed of odd numbers in $[1,n]$. 
Now, let $\tilde{S}(n)=S(n) \setminus (S(\frac{1}{3}n,n) \cup S^1(n))$. Sapozhenko shows that $\tilde{S}(n)$ has an almost correct system of containers. Hence, there exists a subfamily of $\tilde{S}(n)$, $\tilde{S}_B(n)$, with a correct systems of containers. Thus, $|\tilde{S}(n)\setminus \tilde{S}_B (n)|=o(2^{\frac{1}{2}n})$. Sapozhenko proves that $|\tilde{S}_B (n)|=o(2^{\frac{1}{2}n})$, and by consequence $|\tilde{S} (n)|=o(2^{\frac{1}{2}n})$. So, $S(n) \setminus (S(\frac{1}{3}n,n) \cup S^1(n)) = o(2^{\frac{1}{2}n})$ leading to Theorem \ref{thm:Sapozhenko_ce_conjecture}.

One could be forgiven for thinking that given the conjecture by Cameron and \Erdos has been proven (so, it is a theorem now) mathematicians would no longer be interested. However, in 2014 Alon et al. \cite{alon2014refinement} came back to this problem from a different perspective. This time presenting a refinement rather than a proof of the old conjecture. This takes the following form.

\begin{thm}
\label{thm:alon_ce_refinement}
There exists a constant $c>0$ such that for every $n \in \mathbb{Z}^+$ and every $1 \leq m \leq \lceil \frac{1}{2}n \rceil$, the interval $[1,n]$ contains at most $2^{c\frac{n}{m}} \binom{\lceil \frac{1}{2}n\rceil}{m}$ sum-free sets of size $m$.
\end{thm}

The above presents a considerable difference in relation to Theorems \ref{thm:green_ce_conjecture} and \ref{thm:Sapozhenko_ce_conjecture}, as the bound is now specific to the cardinality of the sum-free sets. In order to prove Theorem \ref{thm:alon_ce_refinement} Alon et al. make use of a general bound on the number of independent sets with cardinality $m$ in $3-$uniform hypergraphs (this bound has been proven by the same authors, see  \cite{alon2014counting}). Such sets have a clear relation to sum-free sets. Alon et al. get a hypergraph to encode Schur triples, and eventually proved Theorem \ref{thm:alon_ce_refinement}.

More recently Hancock et al. \cite{hancock2019independent} worked on a slightly different counting problem. How many subsets $S$ of $[1,n]$, sum-free or not, can be partitioned into two sum-free sets (this is sometimes referred to as a 2-wise sum-free set). They put forward the following theorem.

\begin{thm}
\label{thm:n_paritioned}
The number of integer subsets of $[1,n]$ that can be partitioned into two sum-free sets is $\Theta(2^{\frac{4}{5}n+o(n)})$.
\end{thm}

The above means that the number of these subsets grows asymptotically as fast as $2^{\frac{4}{5}n+o(n)}$, which is a tighter bound than if using $\mathcal{O}$. The proof for the above, very much like that for Theorem \ref{thm:Sapozhenko_ce_conjecture}, makes use of advancements in container theory. Perhaps, the most important sum-free set related outcome of Hancock et al. is the conjecture that the number of sets in $[1,n]$ that can be partitioned into two sum-free sets is $\Theta(2^{\frac{4}{5}n})$. This conjecture was finally proven by Tran \cite{tran2018structure} with the following theorem.

\begin{thm}
\label{thm:2_wise_sets_tran}
The number of 2-wise sum-free sets of $[1,n]$ is $\mathcal{O}(2^{\frac{4}{5}n})$
\end{thm}

Clearly, $2^{\frac{4}{5}n} > 2^{\frac{1}{2}n}$ (this inequality is strictly greater because $n \geq 1$). The fact there are more 2-wise sum-free sets in $[1,n]$ than sum-free sets is not particularly surprising. We have that: 
\begin{enumerate}
    \item Every subset of a sum-free set is sum-free, so all sum-free sets with cardinality two or greater can be partitioned into two sum-free sets.
    \item Some sets that are not sum-free can be partitioned into two sum-free sets.
\end{enumerate}

Hence, point 1 gives us $n$ sets that are not 2-wise sum-free in $[1,n]$ (the singletons). However, point 2 gives us more than $n$ 2-wise sum-free sets. There are many examples we could give to show this is indeed the case. For instance, we can join any singleton $S$ containing an even number to any set composed of odd numbers, and partition this union into two sum-free sets (one with the even number and the other with the odd numbers). For any singleton $S$ containing an odd number we can join $S$ to any set composed of odd numbers that is disjoint to $S$, this union would also be 2-wise sum-free.

The proof of Theorem \ref{thm:2_wise_sets_tran} also makes use of the method of containers. The fact this method has been used so often to count sum-free sets, leading to new bounds (see for instance \cite{sapozhenko2003cameron,sapozhenko2008cameron,hancock2019independent,tran2018structure}), is unsurprising. Tan explains this is a powerful tool to deal with combinatorial problems. In essence this shows that independent sets tend to cluster together in many types of hypergraphs. The advantage is that this allows counting the independent sets one cluster at a time.

\begin{figure}
     \centering
     \begin{subfigure}[b]{0.47\textwidth}
         \centering
         \includegraphics[width=\textwidth]{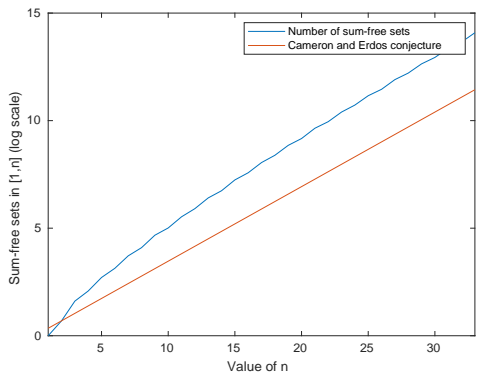}
         \caption{The quantity of sum-free sets for each value of $n$ (blue), and the bound conjectured by Cameron and \Erdos (orange).}
         \label{fig:Number_of_integer_sum_free_against_ce_conjecture}
     \end{subfigure}
     \hfill
     \begin{subfigure}[b]{0.47\textwidth}
         \centering
         \includegraphics[width=\textwidth]{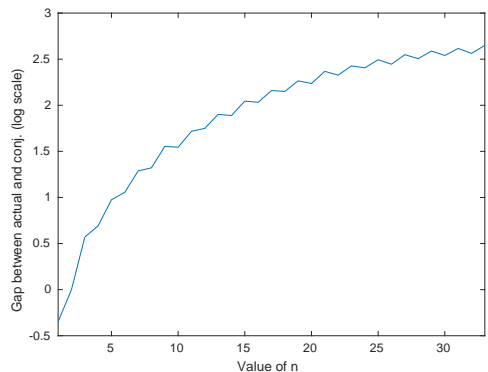}
         \caption{The difference between the number of sum-free sets for each $n$ and the bound conjectured by Cameron and \Erdos.}
         \label{fig:difference_sum_free_against_ce_conjecture}
     \end{subfigure}
     \hfill
        \caption{A comparison, using the logarithmic scale, between the number of sum-free sets for a particular $n$ (for $n=1, 2, \ldots 33$) and the bound conjectured by Cameron and \Erdos.}
        \label{fig:comparing_sum_free_against_ce_conjecture}
\end{figure}

We conclude this section by empirically showing how close the bound conjectured by Cameron and \Erdos actually is in a small sample. Using Algorithm \ref{alg:all_sum_free} we counted all sum-free sets for $n=1,2, \ldots, 33$. Figure \ref{fig:Number_of_integer_sum_free_against_ce_conjecture} shows two curves in the logarithmic scale, the number of sum-free sets for a particular $n$ (in blue) and the bound conjectured by Cameron and \Erdos (in orange). These two curves seem to become parallel as $n$ grows, which is what one would expect. Figure \ref{fig:difference_sum_free_against_ce_conjecture} shows the difference between the two curves in Figure \ref{fig:Number_of_integer_sum_free_against_ce_conjecture}. The ``spiky'' shape in Figure \ref{fig:difference_sum_free_against_ce_conjecture} is also expected as the convergence depends on the parity of $n$ (see Theorem \ref{thm:green_ce_conjecture}).

\subsection{Counting sum-free subsets of abelian groups}
\label{sec:counting_sum_free_in_general}

In Section \ref{sec:ce_conjecture} we considered Cameron and \Erdos proposition that if $S$ is a sum-free integer subset of $[1,n]$, then $|S| \leq \lceil \frac{1}{2}n \rceil$. Kedlaya \cite{kedlaya1998product} shows, with a simple proof, this also applies to finite groups.

\begin{prop}
\label{prop:trivial_bound_sum_free_group}
Let $S$ be the largest sum-free subset of a finite group $G$, then $|S|\leq \frac{1}{2}|G|$.
\end{prop}
\begin{proof}
Let $g \in S$. Given $S$ is sum-free we have that $(S+g) \cap S = \emptyset$. Hence, $2|S|\leq |G| \iff |S| \leq \frac{1}{2}|G|$.
\end{proof}

Recall that given a proper subset $H$ of a finite abelian group $G$, if $g \in G\setminus H$ then the coset $H+g$ is sum-free (see Lemma \ref{lemma:some_coset_is_sum_free}). Using this, Babai and S\'os \cite{babai1985sidon} refined the bound in Proposition \ref{prop:trivial_bound_sum_free_group} by proving the below.

\begin{prop}
\label{prop:coset_bound_sum_free_in_group}
Let $G$ be a group of order $n$, $H$ be a proper subgroup of $G$, and $S$ be the largest sum-free subset of $G$. If $|G:H|=k$ with $k\geq 2$, then $|S|\geq \frac{n}{k}$.
\end{prop}
\begin{proof}
In this case, there are $k$ cosets of cardinality $\frac{n}{k}$. We have that at least one of these cosets follows the format $H+g$ with $g \in G\setminus H$. Lemma \ref{lemma:some_coset_is_sum_free} tells is that any coset under this format is sum-free.
\end{proof}

Using the above propositions we can prove the following Lemma.

\begin{lemma}
\label{lemma:maximum_sum_free_in_index_2}
Let $G$ be an abelian group of order $n$. If $S \subset G$ is sum-free, then $|S|\leq \frac{1}{2}n$ with equality if and only if $S$ is a non-trivial coset of a subgroup of index 2 in $G$.
\end{lemma}
\begin{proof}
Proposition \ref{prop:trivial_bound_sum_free_group} tells us that $|S|\leq \frac{1}{2}n$ does hold. Regarding the equality statement, if $S$ is a non-trivial coset of a subgroup of index 2 in $G$, then Proposition \ref{prop:coset_bound_sum_free_in_group} tells us that $|S| = \frac{1}{2}n$. Let us analyse the converse. Let $S$ be a sum-free subset of $G$ with $|S|=\frac{1}{2}n$ and $g_1,g_2,g_3,g_4 \in S$. Clearly, $g_1+g_2 \in G\setminus S$ and $g_3+g_4 \in G \setminus S$. Notice that $(g_1+g_2)+(g_3+g_4) \in G\setminus S$ as well. Hence, by the finite subgroup criteria $G\setminus S$ is a subgroup of $G$. Let $H$ be this subgroup, we have that $|H| = \frac{1}{2}n$ and $|G:H|=2$. Hence, $S$ is in fact a coset $H+g$ with $g \in G\setminus H$.
\end{proof}

The above can be helpful in finding the largest sum-free subset of a group $G$, which in turn can be helpful in counting the number of sum-free subsets in $G$ (a subset of a sum-free set is sum-free) but let us go in order.

We previously discussed a result by Alon \cite{alon1991independent} that directly applies to any finite abelian group (see Corollary \ref{cor:alon_sum_free_sets_bound}). His corollary implies there are at most $2^{(\frac{1}{2}+o(1))|G|}$ sum-free subsets in any finite abelian group $G$. Lev, Łuczak and Schoen \cite{lev2001sum}, extend the work of Alon by providing a sharp result still. 

\begin{thm}
\label{thm:lev_bound_finite_abelian}
Let $G$ be an abelian group of order $n$. There is an absolute constant $c>0$ such that the number of sum-free sets in $G$ is $(2^{\mathcal{V}(G)}-1)2^{\frac{1}{2}n}+\mathcal{O}(2^{(\frac{1}{2}-c)n})$, where $\mathcal{V}(G)$ is the number of even order components in the canonical decomposition of $G$ into a direct sum of its cyclic groups, and the implicit constant in the $\mathcal{O}-$sign is absolute.
\end{thm}

$G$ is a finite abelian group, and we assume non-trivial, so it can be decomposed to the direct product of cyclic p-groups. This decomposition is of the form $G \cong C_{p_1}^{m_1}, \ldots, C_{p_r}^{m_r}$, where each $p_i$ with $1\leq i \leq r$ is a prime number. The only even prime is 2, so $\mathcal{V}(G)=m_i$ for $p_i=2$. Hence, Theorem \ref{thm:lev_bound_finite_abelian} applies solely to finite abelian groups of even order.

In order to prove Theorem \ref{thm:lev_bound_finite_abelian} the original authors first prove that the number of ``primitive'' sum-free subsets of a finite abelian group $G$ of order $n$ is $\mathcal{O}(2^{(\frac{1}{2}-c)n})$ for an absolute constant $c>0$. A sum-free set $S$ is primitive if there is a homomorphism mapping $S$ to a sum-free $\bar{S} \in G/H$, where $H \subseteq G$. They then combine this result, a lemma that there are at most $2^{(\text{log}_2n)^2}$ subgroups in $G$, and a second lemma stating there are $2^{\mathcal{V}(G)}-1$ subgroups of index two in $G$. Noticing that if $[G:H]=2$, then any sum-free sets in $G$ that are not primitive are all subsets of $G\setminus H$ (this is a consequence of Lemma \ref{lemma:maximum_sum_free_in_index_2}), they eventually prove Theorem \ref{thm:lev_bound_finite_abelian}.

Green and Ruzsa \cite{green2005sum} go further in the analysis of sum-free subsets of finite abelian groups. In their work they first identify the maximum density of a sum-free set.

\begin{defn}
\label{def:density}
Let $S$ be the largest sum-free subset of a finite abelian group $G$, and $\mu(G)$ the density of $S$. That is, $\mu(G) = \frac{|S|}{|G|}$ and $|S| =\mu(G)|G|$. We call $\mu(G)$ the maximum density of a sum-free set in $G$.
\end{defn}

Green and Ruzsa analyse the inequality $\mu(C_m) \geq \frac{1}{m}\lfloor \frac{m+1}{3}\rfloor$, and using Rhemtulla and Street \cite{rhemtulla1970maximal} result that $\mu(C_7^m)=\frac{2}{7}$ for all $m$, they conclude that $\mu(G)\geq \frac{2}{7}$ holds for any finite abelian group. Green and Ruzsa define a function mapping any given finite abelian group to the density of its largest sum-free set, as follows.

\begin{defn}
\label{defn:function_v_green_ruzsa}
Let $v$ be the function $v:G\rightarrow [\frac{1}{2}, \frac{2}{7}]$, where $G$ is any finite abelian group, be defined as follows.
\[ v(n)=\begin{cases} 
      \frac{1}{3}+\frac{1}{3p} &, \text{if }|G| \text{ is divisible by a prime } p\equiv 2 \text{(mod 3), }p\text{ is the lowest such prime}.\\
      \frac{1}{3} &, \text{the above does not apply but } 3 \mid |G|.\\
      \frac{1}{3}-\frac{1}{3m} &, |G| \text{ is divisible only by primes } p\equiv 1\text{ (mod 3), and }m\text{ is the largest}\\
      &\text{order of any element of }G.
   \end{cases}
\]
\end{defn}

The authors then divide sum-free sets in three types (matching the three cases in the definition of $v$, above), and eventually prove that $\mu(G) = v(G)$.

The above gives us the cardinality of the largest sum-free set in any finite abelian group $G$ (something we take further when discussing maximal sum-free sets, in Section \ref{chap:maximal}). Noting that any subset of a sum-free set is sum-free, Green and Ruzsa carry on. Let $f(G)$ be the number of sum-free sets in $G$, then $ f(G) \geq 2^{\mu(G)|G|}$. That is, the number of sum-free sets in $G$ is at least as high as the number of subsets in its largest sum-free set. They then show that nearly all sum-free sets in $G$ are contained in some sum-free subset of maximum cardinality. The general idea that nearly all sum-free sets are in fact subsets of a small number of other sum-free sets is not completely alien. Green \cite{green2004cameron} and Sapozhenki \cite{sapozhenko2008cameron} had previously shown that nearly all integer sum-free subsets of $[1,n]$ are either entirely composed of odd numbers, or contained in the interval $[\lceil\frac{n+1}{3} \rceil,n]$ (see Section \ref{sec:proof_ce_conjecture}). With the above, Green and Ruzsa prove the following theorem.

\begin{thm}
\label{thm:bound_sf_finite_g}
Let $G$ be a finite abelian group, $f(G)$ be the number of sum-free subsets of $G$, and $\mu (G)$ be the density of the largest sum-free subset of $G$. Then, $f(G)=2^{(1+o(1))\mu(G)|G|}$.
\end{thm}

Notice that calculating the above is relatively easy given the existence of function $v$ (see Definition \ref{defn:function_v_green_ruzsa}). 

It can also be interesting to count sum-free subsets carrying a particular property. With this in mind Alon et al. \cite{alon2014counting} extend the work of Green and Ruzsa by counting sum-free sets with a particular cardinality $m$, as we can see in their theorem below. 

\begin{thm}
\label{thm:alon_sum_free_card_m}
Let $G$ be an abelian group of order $n$. Then, 
\begin{align*}
    f(G,m) = (\#\{\text{elements of }G\text{ of order }2\} +o(1)) \binom{n/2}{m}.
\end{align*}
\end{thm}

In the above $f(G,m)$ is the number of sum-free subsets of $G$ with cardinality $m$. Alon et al. proof of Theorem \ref{thm:alon_sum_free_card_m} is rather long and involved. They are able to determine an upper-bound for independent sets of cardinality $m$ on a particular type of graph. They then partition sum-free subsets of $G$ in two sets based on the intersection between these sets and the subgroups of $G$ with index 2. Modelling the sets using graphs, Alon et al. are able to use the found upper-bound to eventually prove Theorem \ref{thm:alon_sum_free_card_m}.


\section{Maximal sum-free sets}
\label{chap:maximal}

In the previous section we saw that Green and Ruzsa \cite{green2005sum} proved that nearly all sum-free subsets of a finite abelian group $G$ are contained in a sum-free subset of $G$ of maximum cardinality (see the discussion of Theorem \ref{thm:bound_sf_finite_g}). This sum-free subset of maximum cardinality belongs to a well-known class of sum-free sets called maximal sum-free.

\begin{defn}
Let $G$ be a group and $S\subset G$ be sum-free. We call $S$ maximal sum-free if there is no sum-free $A\subset G$ such that $S\subset A$. That is, $S$ is not properly contained in any other sum-free subset of $G$. 
\end{defn}

\subsection{Maximal sum-free subsets of positive integers}

Cameron and \Erdos \cite{cameron1999notes} investigate maximal sum-free integer subsets of $[1,n]$. They explain they expect the number of maximal sum-free sets, $f_{\text{max}}(n)$, to be substantially lower than that of sum-free sets. This expectation seems to be supported by Theorem \ref{thm:Sapozhenko_ce_conjecture} and the proof of Theorem \ref{thm:green_ce_conjecture}, which imply the almost all sum-free integer subsets of $[1,n]$ are either composed entirely of odd numbers (i.e. subsets of a maximal sum-free set of all odd numbers) or are contained in the interval $[\frac{1}{3}n,n]$. The latter is not sum-free but it contains the maximal sum-free sets 
$[\frac{1}{2}(n+1),n]$ (for an odd $n$), and $[\frac{1}{2}n,n-1]$ as well as $[\frac{1}{2}n+1,n]$ (both of them for an even $n$), identified in their previous work (see \cite{cameron1990number}, and the discussion in Section \ref{sec:ce_conjecture}). Cameron and \Erdos then prove the following.

\begin{thm}
The number of maximal sum-free integer subsets of $[1,n]$ is equal or greater than $2^{\lfloor \frac{n}{4}\rfloor}$.
\end{thm}

Interestingly, their proof of the above does not require the identification of all maximal sum-free sets of $[1,n]$. This is probably easier to explain with an example (an interested reader can find the proof itself in \cite{cameron1999notes}).

\begin{example}
Le us analyse $[1,n]$ with $n=6$. Let $m=n$ or $m=n-1$, whichever leads to an even $m$. Hence, $m=n=6$. Let us now build sum-free sets consisting of $m$ and only one component of each pair $(x,m-x)$ for any odd $x<\frac{m}{2}$. The only suitable value of $x$ is one, leading to the pair $(1,5)$. With this we get the sum-free sets $\{1,6\}$ and $\{5,6\}$. Notice that these sets are not necessarily maximal, however, no further odd number less than $m$ can be added to any of them. For instance, if you take $\{1,6\}$ and add to it either $3$ or $5$ this set will cease to be sum-free. Hence, both $\{1,6\}$ and $\{5,6\}$ are subsets of different maximal sum-free sets. In this example we found the maximal sum-free sets (among others): $\{4,5,6\}$, $\{2,5,6\}$, and $\{1,4,6\}$. These show that (i) the sets $\{1,6\}$ and $\{5,6\}$ are indeed subsets of different maximal sum-free sets; (ii) a set may in fact be a subset of more than one maximal sum-free set.
\end{example}

The problem of finding an upper bound for $f_{\text{max}}(n)$ has attracted some attention. \L uczak and Schoen \cite{luczak2001number} devised a probabilistic proof showing that $f_{\text{max}}(n) \leq 2^{\frac{1}{2}n-2^{-28}n}$. Notice that every maximal sum-free set is of course sum-free, and the latter grows at $\mathcal{O}(2^{\frac{1}{2}n})$ in $[1,n]$ (see Section \ref{sec:proof_ce_conjecture}). Hence, the result above by \L uczak and Schoen is an improvement but intuition would indicate there is room for more. Wolfovitz \cite{wolfovitz2009bounds} took this further and proved that $f_{\text{max}}(n) \leq 2^{\frac{3}{8}n+o(n)}$. His proof made use of the family of sum-free sets $\mathcal{F}$ devised by Green \cite{green2004cameron} when proving Theorem \ref{thm:green_ce_conjecture} (for a description, see Section \ref{sec:proof_ce_conjecture}), and another family of sum-free sets whose union is maximal sum-free. Using these he eventually reaches the stated bound.

Balogh et al. \cite{balogh2015number} make a considerable improvement on the above with the following theorem.

\begin{thm}
\label{thm:maximal_upperbound_intefer_sum_free}
There are at most $2^{(\frac{1}{4}+o(1))n}$ maximal sum-free integer subsets of $[1,n]$.
\end{thm}
The proof of the above makes use of containers and the family of sum-free sets $\mathcal{F}$ (see Section \ref{sec:proof_ce_conjecture}). Container theory is in fact a quite popular approach to count sum-free sets, maximal or not. We have seen this being used by Sapozhenko \cite{sapozhenko2008cameron}, Hancock et al. \cite{hancock2019independent}, and Tran \cite{tran2018structure} when proving Theorems \ref{thm:Sapozhenko_ce_conjecture}, \ref{thm:n_paritioned}, and \ref{thm:2_wise_sets_tran}, respectively. Each element of $\mathcal{F}$ is a container. Balogh et al. then prove that each one of these containers has at most $2^{\frac{1}{4}n+o(n)}$ maximal sum-free sets, leading to Theorem \ref{thm:maximal_upperbound_intefer_sum_free}.

Recently, Balogh et al. \cite{balogh2018sharp} took the above even further by providing an exact answer, instead of a bound, to the number of maximal sum-free sets. They prove the theorem below using containers and the family of sum-free sets $\mathcal{F}$ (see Section \ref{sec:proof_ce_conjecture}). However, they were forced to make improvements in container theory in order to avoid over-counting the number of maximal sum-free sets thanks to an error term in the original theory.

\begin{thm}
\label{thm:number_maximal_sum_free_in_integers}
For each $1\leq i \leq 4$, there is a constant $c_i$ such that given any $n \equiv i$ (mod $4$), $[1,n]$ contains $(c_i+o(1))2^{\frac{1}{4}n}$ maximal sum-free sets.
\end{thm}

The proof of the above also sheds some light on the structure of maximal sum-free subsets of $[1,n]$. For instance, if $S \in [1,n]$ contains $o(n^2)$ Schur triples (such triple consists of $x,y,z \in S$ such that $x+y=z$), then $S$ follows under one of three categories:
\begin{enumerate}
    \item $|S| \leq 0.47n$
    \item $|S| = (\frac{1}{2}-\gamma)n$ and $S = A \cup B$ where $|A|=o(n)$ and $B\subseteq [(\frac{1}{2}-\gamma)n,n]$ is sum-free.
    \item $|S| = (\frac{1}{2}-\gamma)n$ and $S = A \cup B$ where $|A|=o(n)$ and $B$ contains solely odd numbers.
\end{enumerate}

Notice that $S$ is not necessarily sum-free, however, a family of such sets would be nearly sum-free. With this, Balogh et al. were able to describe general properties of maximal sum-free sets. Here, we state those we found the most interesting. There are only a small number of maximal sum-free subsets of $[1,n]$ in category 1. This seems well-aligned with intuition. The cardinality of every sum-free set of maximum cardinality is $\lceil \frac{1}{2}n \rceil$, and these are maximal sum-free sets as well. One would then expect that the majority of maximal sum-free sets, that are not of maximum cardinality, to be not much lower than $\lceil \frac{1}{2}n\rceil$. Regarding category 2, Balogh et al. show that the number of maximal sum-free sets under this category with two or more even numbers amounts to only $o(2^{\frac{1}{4}n})$. This also follows intuition. Such set would have a cardinality higher than $0.47n$ (otherwise it would be category 1), and the sum of two even numbers is even. The more even numbers a sum-free set has, the hardest it is to achieve such high cardinality. Finally, by analysing the sets in categories 2 and 3 they find that the number of maximal sum-free subsets of $[1,n]$ containing at most one even number is $\mathcal{O}(2^{\frac{1}{4}n})$. This is rather close to the bound established in Theorem \ref{thm:maximal_upperbound_intefer_sum_free}. Hence, the majority of maximal sum-free sets falls under these categories and follow this configuration.

\subsection{Maximal sum-free subsets of abelian groups}

In this section we explore maximal sum-free subsets of finite abelian groups. Here, $f_{\text{max}}(G)$ represents the number of maximal sum-free sets in a group $G$.

\begin{thm}
\label{thm:maximal_sum_free_previous}
Let $G$ be a group of order $n$, $f_{\textup{max}}(G) \leq 2^{0.406n+o(n)}$.
\end{thm}

Wolfovitz \cite{wolfovitz2009bounds} proved the above by first identifying a sum-free $S \subset G$, and constructing a Cayley graph. This graph is $k$-regular (i.e. each vertex has degree $k$) with $k=|S\cup (-S)|$. Wolfovitz carries on by stating there are at most $2^{o(n)}$ maximal sum-free subsets of $G$ with cardinality less than $\lceil \sqrt{n} \rceil$ and concludes that to prove Theorem \ref{thm:maximal_sum_free_previous} it suffices to count those with cardinality of at least $\lceil \sqrt{n} \rceil$. Each maximal sum-free subset is an independent set in at least one Caylay graph. Since there are $\binom{n}{\lceil \sqrt{n} \rceil}=2^{o(n)}$ such graphs, it is enough to fix a sum-free $S \subset G$ of size $\lceil \sqrt{n} \rceil$ and upper bound the number of independent sets that correspond to the maximal sum-free sets. This eventually leads to Theorem \ref{thm:maximal_sum_free_previous}.

Balogh et al. \cite{balogh2018sharp} argues that Theorem \ref{thm:maximal_sum_free_previous} only shows that $f_{\max}(G)$ is exponentially smaller than $f(G)$ if $G$ is of even order (as its largest sum-free subset will have cardinality $\frac{1}{2}n$). They then produce a proof for the improvement below.

\begin{thm}
\label{thm:maximal_sf_bound_upper}
Let $G$ be an abelian group of order $n$, and $S$ be the largest sum-free subset of $G$. Then, $f_{\textup{max}}(G) \leq 3^{\frac{|S|}{3}+o(n)}$.
\end{thm}

To prove the above, Balogh et al. constructs a family $\mathcal{F}$ of containers that is slightly different than that devised by Green to prove Theorem \ref{thm:green_ce_conjecture}. Recall that $\mu(G) = \frac{|S|}{|G|}$, where $S$ is the largest sum-free subset of a group $G$ (see Definition \ref{def:density}). In this new version we have that for any $F \in \mathcal{F}$, $F=B\cup C$, where $B$ is sum-free in $G$ and $|B|\leq \mu(G)|G|$, and $C \subset G$ with $|C|=o(n)$. The other previous requirements of $|\mathcal{F}|=2^{o(n)}$ and that any sum-free subset of $G$ is contained in at least one element of $\mathcal{F}$ remain. A sum-free $B\subset G$ cannot have a cardinality that is higher than that of the largest sum-free subset of $G$. 

Balogh et al. then fixes $F \in \mathcal{F}$ and states that every maximal sum-free subset of $G$ contained in $F$ can be formed by picking a sum-free $S\subseteq C$ (at most $2^{o(n)}$ choices, by the definition above) and extending it with $B$. Regarding the latter, Balogh et al. state that its number of maximal independent sets in a graph is at most $3^{\frac{|B|}{3}} \leq 3^{\frac{\mu(G)|G|}{3}}$. These two results eventually lead to Theorem \ref{thm:maximal_sf_bound_upper}.

If a sum-free set $S$ is of maximum cardinality, then $S$ is maximal sum-free. Such $S$ is of particular importance because nearly all sum-free subsets of the same group are contained in a set like $S$ (see the discussion of Theorem \ref{thm:bound_sf_finite_g}). Also, in the case of subsets of $[1,n]$ Balogh et al. \cite{balogh2018sharp} showed that only a small number of maximal sum-free sets have cardinality less or equal to $0.47n$ (see our discussion of Theorem \ref{thm:number_maximal_sum_free_in_integers}). Hence, it is interesting to study the exceptions. That is, groups with small maximal sum-free sets.

Giudici and Hart \cite{giudici2009small} explore the world of small maximal sum-free sets. Their work proves various interesting theorems. Here, we present some of their main findings.

\begin{thm}
\label{thm:maximal_of_cardinality_one}
Let $G$ be a group and $S\subset G$ be a maximal sum-free set of cardinality one. Then, $G \cong C_2, C_3, C_4$ or $Q_8$, and the element in $S$ is of prime order in $G$.
\end{thm}

Giudici and Hart prove the above by analysing $S$, a maximal sum-free set containing a single element $g$. They notice that if $g$ does not have order 2, then $S \cap -S=\emptyset$. This is clear as if $o(g)=2$ we would have $-g=g$, by consequence $S=-S=\{g\}$ and $S\cap -S=\{g\}$. Using the $S \cap -S=\emptyset$ result and a previous corollary in their work \cite[Corollary 3.10]{giudici2009small} stating that for this case $|G|=4\times |S|^2+1$, they reach $|G|\leq 5$. Then, after testing all possibilities they conclude $G \cong C_3$. Afterwards they examine the case $o(g)=2$. In this case they find that every $x \in G\setminus \langle g \rangle$ has order 4 and $\langle g \rangle$ is the only subgroup of $G$ with order 2. They then conclude $|G| \in \{2,4,8\}$ and $G\cong C_2, C_4$ or $Q_8$. 

Giudici and Hart also prove that only certain groups have a sum-free subset of cardinality 2. They provide a list containing 11 such groups and their respective sum-free sets (note that a group may have more than one sum-free set with cardinality 2). $C_4$ is one of the groups in the list, as one would expect (given that $\frac{1}{2}|C_4|=2$) but it is very interesting to see there are various other groups listed, one with order 16. Finally, Giudici and Hart conjecture that any group $G$ with $|G|>24$ does not contain a maximal sum-free set of size 3. This conjecture was later proved by Anabanti and Hart \cite{anabanti2015locally}.

\section{Conclusion}
\label{chapter:conclusion}


In this paper our main aim is to present and discuss some interesting properties of sum-free subsets of abelian groups, to a broader readership. More specifically, we set out to discuss the answers to questions like those raised by Alon et al. \cite{alon2014counting} and Green and Ruzsa \cite{green2005sum}: how many sum-free subsets are there in a given abelian group $G$? what are its sum-free subsets of maximum cardinality? what is the maximum cardinality of these sum-free subsets? what does a typical sum-free subset of $G$ looks like?

In order to address the questions above, and others, we first asked ourselves: given an abelian group $G$ how can we identify its sum-free subsets? Section \ref{chap:Identifying_sum_free_sets} addresses this by presenting different methods, computational and algebraic, that can be employed to identify such sum-free sets. In Section \ref{chapter:counting_sum_free} we discussed a considerable amount of work done on counting sum-free sets. We also identified that Cameron and \Erdos\kern-0.65ex' list of integer sum-free subsets of $[1,n]$ with cardinality $\lceil \frac{n}{2}\rceil$ is incomplete (we found exceptions for three values of $n$ and these are the only exceptions there are, see Section \ref{sec:ce_conjecture}). Section \ref{chapter:counting_sum_free} also explained that nearly all sum-free subsets of $G$ are contained in some sum-free set of maximum cardinality (see the discussion of Theorem \ref{thm:bound_sf_finite_g}). This clearly showed the importance of sum-free sets of maximum cardinality, giving us reason to focus Section \ref{chap:maximal} on maximal sum-free sets. Section \ref{chap:maximal} discussed proofs regarding the bounds of maximal sum-free sets. It also sheds some light on the structure of maximal sum-free sets (see Theorem \ref{thm:number_maximal_sum_free_in_integers}). Within these sections, we have covered what we think are the main areas pertinent to sum-free subsets of abelian groups.




 \textbf{Data availability statement} Data sharing not applicable to this article as no data sets were generated or analysed during the current study.
\bibliographystyle{ieeetr}
\bibliography{references}

\begin{thebibliography}{10}

\bibitem{schur1916}
I.~Schur, ``Uber die kongruenz {$x^m + y^m \equiv z^m ($ mod $ p)$},'' {\em
  ber. Deutsch. Mat. Verein.}, no.~25, pp.~114--117, 1916.

\bibitem{cameron1987portrait}
P.~J. Cameron, {\em Portrait of a typical sum-free set}, vol.~123.
\newblock Cambridge Univ. Press, 1987.

\bibitem{cameron1990number}
P.~J. Cameron and P.~\Erdos, ``On the number of sets of integers with various
  properties,'' {\em Number Theory (RA Mollin, ed.)}, pp.~61--79, 1990.

\bibitem{alon2014counting}
N.~Alon, J.~Balogh, R.~Morris, and W.~Samotij, ``Counting sum-free sets in
  abelian groups,'' {\em Israel Journal of mathematics}, vol.~199, no.~1,
  pp.~309--344, 2014.

\bibitem{green2005sum}
B.~Green and I.~Z. Ruzsa, ``Sum-free sets in abelian groups,'' {\em Israel
  Journal of Mathematics}, vol.~147, no.~1, pp.~157--188, 2005.

\bibitem{erdijs1965extremal}
P.~\Erdos, ``Extremal problems in number theory,'' in {\em Proc. Sympos. Pure
  Math}, pp.~181--189, 1965.

\bibitem{alon1990sum}
N.~Alon and D.~J. Kleitman, ``Sum-free subsets,'' {\em A tribute to Paul
  \Erdos}, pp.~13--26, 1990.

\bibitem{rhemtulla1971maximal}
A.~Rhemtulla and A.~P. Street, ``Maximal sum-free sets in elementary abelian
  p-groups,'' {\em Canadian Mathematical Bulletin}, vol.~14, no.~1, pp.~73--80,
  1971.

\bibitem{cameron_blog}
P.~J. Cameron, ``{Clive Sinclair}.''
  {http://cameroncounts.wordpress.com/2021/09/18/clive-sinclair/}.
\newblock Accessed on: 2022-08-15.

\bibitem{kolountzakis1994selection}
M.~N. Kolountzakis, ``Selection of a large sum-free subset in polynomial
  time,'' {\em Information processing letters}, vol.~49, no.~5, pp.~255--256,
  1994.

\bibitem{tran2018structure}
T.~Tran, ``On the structure of large sum-free sets of integers,'' {\em Israel
  Journal of Mathematics}, vol.~228, no.~1, pp.~249--292, 2018.

\bibitem{freiman1991structure}
G.~A. Freiman, ``On the structure and the number of sum-free sets,'' {\em
  Journ{\'e}es Arithm{\'e}tiques}, pp.~195--201, 1991.

\bibitem{calkin1990number}
N.~J. Calkin, ``On the number of sum-free sets,'' {\em Bulletin of the London
  Mathematical Society}, vol.~22, no.~2, pp.~141--144, 1990.

\bibitem{szemeredi1975sets}
E.~Szemer{\'e}di, ``On sets of integers containing no k elements in arithmetic
  progression,'' {\em Acta Arith}, pp.~503--505, 1974.

\bibitem{alon1991independent}
N.~Alon, ``Independent sets in regular graphs and sum-free subsets of finite
  groups,'' {\em Israel Journal of Mathematics}, vol.~73, no.~2, pp.~247--256,
  1991.

\bibitem{green2004cameron}
B.~Green, ``The {C}ameron--{E}rd{\H{o}}s conjecture,'' {\em Bulletin of the
  London Mathematical Society}, vol.~36, no.~6, pp.~769--778, 2004.

\bibitem{sapozhenko2003cameron}
A.~Sapozhenko, ``The {C}ameron-{E}rd{\H{o}}s conjecture,'' in {\em Doklady.
  Mathematics}, vol.~68, pp.~438--441, 2003.

\bibitem{sapozhenko2008cameron}
A.~A. Sapozhenko, ``The {C}ameron--{E}rd{\H{o}}s conjecture,'' {\em Discrete
  mathematics}, vol.~308, no.~19, pp.~4361--4369, 2008.

\bibitem{alon2014refinement}
N.~Alon, J.~Balogh, R.~Morris, and W.~Samotij, ``A refinement of the
  {C}ameron--{E}rd{\H{o}}s conjecture,'' {\em Proceedings of the London
  Mathematical Society}, vol.~108, no.~1, pp.~44--72, 2014.

\bibitem{hancock2019independent}
R.~Hancock, K.~Staden, and A.~Treglown, ``Independent sets in hypergraphs and
  ramsey properties of graphs and the integers,'' {\em SIAM Journal on Discrete
  Mathematics}, vol.~33, no.~1, pp.~153--188, 2019.

\bibitem{kedlaya1998product}
K.~S. Kedlaya, ``Product-free subsets of groups,'' {\em The American
  mathematical monthly}, vol.~105, no.~10, pp.~900--906, 1998.

\bibitem{babai1985sidon}
L.~Babai and V.~T. S{\'o}s, ``Sidon sets in groups and induced subgraphs of
  cayley graphs,'' {\em European Journal of Combinatorics}, vol.~6, no.~2,
  pp.~101--114, 1985.

\bibitem{lev2001sum}
V.~F. Lev, T.~Łuczak, and T.~Schoen, ``Sum-free sets in abelian groups,'' {\em
  Israel Journal of Mathematics}, vol.~125, pp.~347--368, 2001.

\bibitem{rhemtulla1970maximal}
A.~Rhemtulla and A.~P. Street, ``Maximal sum-free sets in finite abelian
  groups,'' {\em Bulletin of the Australian Mathematical Society}, vol.~2,
  no.~3, pp.~289--297, 1970.

\bibitem{cameron1999notes}
P.~J. Cameron and P.~Erd{\H{o}}s, ``Notes on sum-free and related sets,'' {\em
  Combinatorics, Probability and Computing}, vol.~8, no.~1-2, pp.~95--107,
  1999.

\bibitem{luczak2001number}
T.~{\L}uczak and T.~Schoen, ``On the number of maximal sum-free sets,'' {\em
  Proceedings of the American Mathematical Society}, pp.~2205--2207, 2001.

\bibitem{wolfovitz2009bounds}
G.~Wolfovitz, ``Bounds on the number of maximal sum-free sets,'' {\em European
  Journal of Combinatorics}, vol.~30, no.~7, pp.~1718--1723, 2009.

\bibitem{balogh2015number}
J.~Balogh, H.~Liu, M.~Sharifzadeh, and A.~Treglown, ``The number of maximal
  sum-free subsets of integers,'' {\em Proceedings of the American Mathematical
  Society}, vol.~143, no.~11, pp.~4713--4721, 2015.

\bibitem{balogh2018sharp}
J.~Balogh, H.~Liu, M.~Sharifzadeh, and A.~Treglown, ``Sharp bound on the number
  of maximal sum-free subsets of integers,'' {\em Journal of the European
  Mathematical Society}, vol.~20, no.~8, pp.~1885--1911, 2018.

\bibitem{giudici2009small}
M.~Giudici and S.~Hart, ``Small maximal sum-free sets,'' {\em the electronic
  journal of combinatorics}, pp.~1--17 (R59), 2009.

\bibitem{anabanti2015locally}
C.~S. Anabanti and S.~Hart, ``Locally maximal product-free sets of size 3,''
  {\em arXiv preprint arXiv:1503.06509}, 2015.

\end{thebibliography}



    

\end{document}